\documentclass[amssymb, 11pt]{amsart}
\usepackage{latexsym}
\usepackage{young}

\newlength{\standardunitlength}
\newcommand{\bea}{\begin{eqnarray}}
\newcommand{\ena}{\end{eqnarray}}
\newcommand{\beas}{\begin{eqnarray*}}
\newcommand{\enas}{\end{eqnarray*}}

\newcommand{\ignore}[1]{}

\newcommand{\SL}{\operatorname{SL}} 
\newcommand{\PSL}{\operatorname{PSL}} 
\newcommand{\PGL}{\operatorname{PGL}} 
\newcommand{\Sp}{\operatorname{Sp}}
\newcommand{\SO}{\operatorname{SO}}

\newcommand{\alt}{\operatorname{Alt}}
\newcommand{\SU}{\operatorname{SU}}

\newtheorem{prop}{Proposition}[section]

\newtheorem{cor}[prop]{Corollary}
\newtheorem{theorem}[prop]{Theorem}

\begin{document}

\title [Involutions in Simple Groups] {Commuting Involutions in Finite Simple Groups}

\author{Robert Guralnick}
\address{Department of Mathematics\\University of Southern California \\ Los Angeles, CA, 90089-2532, USA}
\email{guralnic@usc.edu}

\author{Geoffrey R. Robinson}
\address{Department of Mathematics, King's College, 
Aberdeen,  AB24 3FX, UK}
\email{g.r.robinson@abdn.ac.uk}

\date{\today}

\thanks {Guralnick was partially supported by NSF grant DMS 1901595.}

\begin{abstract}     We prove that if $G$ is a finite simple group and $x, y \in G$ are involutions, then
$|x^G \cap C_G(y)| \rightarrow \infty$ as $|G| \rightarrow \infty$.  This extends results of Guralnick-Robinson
and Skresanov.   We also prove a related result about $C_{G}(t)/O(C_G(t))$ that does not require the classification of finite simple groups.  
\end{abstract}

\maketitle

\section{Introduction}

Let $G$ be a finite simple group.  In \cite{GR}, it was proved that given any pair of involutions $x,y$  in $G$,
$x$ commutes with some conjugate of $y$.  Indeed, it was shown that that there an elementary abelian
subgroup $A$ of $G$ that intersects every conjugacy class of involutions.   One can state this in terms 
of graphs.  Let $X$ be the graph whose vertices are conjugacy classes of involutions of $G$  with two distinct vertices
adjacent if there are commuting elements in the two classes.   The result implies this is a complete graph.   This is not
true for elements of order $p$ for an  odd prime $p$ or for involutions in almost simple groups. 
See \cite{Ca} for  more about such graphs.

Recently,  Skresanov \cite{Sk}  used this to prove that if $x \in G$ is an involution with $G$ a  finite simple group,  then
if the number of involutions centralized by $x$ has size at most $n$ then $|G| \le h(n)$ for some function $h$. 
Moreover, this implies that if $H$ is a  locally finite simple group and $x \in H$ is involution, then $x$ commutes with infinitely
many involutions.  

We generalize this and show that any involution $x$ in a finite simple group must commute with many involutions in every conjugacy class of
involutions.  More precisely, we prove: 

\begin{theorem} \label{t:main}  Let $G$ be a finite simple group.  There exists a function $f$ on the natural numbers
 so that if $x, y \in G$ are involutions and
$|x^G  \cap C_G(y)| \le n$, then $|G| \le f(n)$.
\end{theorem} 

A trivial corollary of this is a result for algebraic groups (which could be proved directly by the same methods).

\begin{cor} \label{c:algebraic}  Let $G$ be a simple  adjoint algebraic group over the algebraically closed field $k$.  If $x, y \in G$ are involutions, then
$\dim x^G \cap C_G(y) > 0$.
\end{cor}

In fact, one can show that  $\dim x^G \cap C_G(y)$ grows linearly in the rank of $G$.   On the other hand, if $G = \PGL_2(k)$, then there is a unique class of involutions
and $\dim x^G \cap C_G(x) = 1$ (and moreover in this case if the characteristic of  $k$ has characteristic not $2$, then $x^G \cap C_G(x)^0 = \{ x \}$).  

Our next result is a variation on the theme for all finite groups with a weaker result but with a proof that does not require
the classification of finite simple groups.  If $H$ is a group, let $\mathcal{I}(H)$ denote the set of involutions in $H$. 

\begin{theorem} \label{t:main2}    Let $n$ be a  positive integer and let $\epsilon > 0$.   There are only finitely many simple groups $G$ with an involution
$t$ such that
$$[C_G(t):O(C_G(t)]=n,$$
and
$$|G| >(1 + \epsilon) |C_G(t)|^2|\mathcal{I}(C_G(t))|.
$$
\end{theorem}

By contrast, for any integer $a >1$, the groups  ${\rm PSL}(2,q)$ with \\
$q \equiv 1+2^{a}$ (mod $2^{a+1}$) provide infinitely many examples of finite simple groups $G$ containing an involution $t$ with  $$[C_{G}(t):O_{2^{\prime}}(C_{G}(t))] = 2^{a}$$ and $$|G| > |C_{G}(t)|^{2}|\mathcal{I}(C_{G}(t))|.$$

\section{Proof of Theorem  \ref{t:main} }

Let $G$ be a finite simple groups with $x, y \in G$ in involutions such that $|x^G \cap C_G(y)| \le n$. 
Since this is an asymptotic result, we can ignore the sporadic groups.    

\subsection{Alternating Groups}   

Let $G=A_m$.   Let $Y_1$ be the set of points moved by $y$ and $Y_2$ the set of fixed points and set $m_i = |Y_i|$.   If $m_1 \ge n/2$,  
we can choose a conjugate $x'$ of $x$ so that $x'$ is nontrivial on $Y_1$ and so that $\langle x', y \rangle$ has an orbit of size $4$.
Working in $H = \alt(Y_1)$, we see that the  $C_G(y)$ orbit containing $x'$ grows with $m$.   If $m_2 \ge n$, we can choose $x'$ so that it acts nontrivially
on $Y_2$ and so conjugating by $\alt(Y_2) \le C_G(y)$ gives a large orbit on $x^G$ (one that grows with $m$).  
Thus, the result holds in this case. 

\subsection{Classical Groups in Odd Characteristic}

First consider the case of $\PSL_2(q)$ with $q$.  Then $G$ contains a unique class of involutions and the centralizer of an involution
is a dihedral group of order $q \pm 1$ and the number of involutions in the centralizer is $1 + (q \pm 1)/2$ and grows with $q$. 

In the remaining cases, we  work in the groups $\SL_n(q), n \ge 3, \SU_{n}(q), n \ge 3, \Sp_n(q), n \ge 4$ and $\SO^{\epsilon}_n(q), n \ge 7$.  

If the involution lifts to an involution in the linear group (or more general has eigenvalues in the field of definition for the natural module).  then we argue
just as for the alternating groups -- we decompose the natural module $V$ as direct sum of the eigenspaces of $y$
and note that the derived group  $D$ of $C_G(y)$ is a direct product of two groups of the same type and at least one
of them is large.  We choose $x'$ a conjugate of $x$ so that $x'$ commutes with $y$ but does not commute with $D$
Then $|(x')^D|$ grows as $q$ or $n$ (or both) grow and the result follows.

The other possibility is that $y^2$ is a scalar in the linear group but the eigenvalues of $y$ are not in the base
field.  This case the derived group $D$ of $C_G(y)$ is  an extension field group.
Again, we choose $x' \in C_G(y)$ that does not centralizes $D$ and so $|(x')^D|$ grows with $n$ or $q$.

\subsection{Classical Groups in Characteristic $2$} 

The proof is quite similar to the proof given of the existence of broad subgroups given in \cite{GR}. We use the classification
of involutions (see  (\cite[]{LS}, \cite{AS}, \cite{He}).   We note that every involution has a fixed space of at least half the dimension f
or any module in characteristic $2$.

First we consider the cases where $G$ is not an orthogonal group.  Then every conjugacy class of involutions intersects $Q$, the unipotent radical of
the stabilizer $P$ of 
a maximal isotropic subspace (in the case of $\SL_n$ we take a subspace of dimension $n/2$ or $(n-1)/2$).   
Note that $Q$ is abelian and that if $L$ is the Levi subgroup of the parabolic subgroup $P=QL$,
then every orbit of $L$ grows with the field size and dimension.  Thus  $|x^G \cap C_G(y)| \ge |x^G \cap Q| \rightarrow \infty$
as $q + n \rightarrow \infty$.  

Next consider $G=\SO^{\pm}_{2n}(q)$ with $q$ even and $n \ge 4$.   Let $V$ denote the natural module for $G$.  If $n$ is odd, then the
fixed space of any involution on $V$ has dimension least $n+1$ (since the number of nontrivial Jordan blocks is even)
and so cannot be totally singular and so stabilizes a nondegenerate $1$-space.
The stabilizer of a nondegenerate $1$-space in the full orthogonal group is $C \times \Sp_{2n-2}(q)$  where $C$ has order $2$ and is generated by a transvection.
Thus every conjugacy of involutions intersects $C \times \Sp_{2n-2}(q)$ and its projection on the second factor is nontrivial.  
The result for symplectic groups yields the desired conclusion.  Similarly if 
$\epsilon$ is of $-$ type, the fixed space has dimension at least $n$ and a maximal isotropic subspace has dimension $n-1$ and the same argument
applies.  

The remaining case is when the orthogonal group has plus type and $n$ is even.  From the classification of involutions \cite{AS, He, LS}, it follows that
we can decompose $V = A \perp B$ where $A$ is $4$-dimensional of $+$ type and $x, y$ preserve the decomposition, commute on $A$ and
each preserves a totally singular subspace in the same $\SO_{2n-2}(q)$ orbit.  Now apply the result for $\SO_{2n-2}(q)$ (if $n =4$, then one would need
to check this for $\SO_4(q)$ or alternatively, we may assume that $q$ is increasing and then we see that $x, y \in \SO_4(q) \times SO_4(q)$ and 
the Sylow $2$-subgroups of these subgroups are abelian and have order $q^4$.   Since every involution is conjugate to an element in
this subgroup, the result follows). 

\subsection{Exceptional Groups in Odd Characteristic}  

The classes of involutions in the exceptional groups are all given in \cite{Lu}.   We go through the various cases.  

If $G={^2}G_2(q)$ with $q= 3^{2a+1} > 3$, there is a unique class of involutions and the centralizer is $2 \times \PSL_2(q)$
and the result holds.  If $G={^3}D_4(q)$, there is a unique class of involutions with centralizer $\SL_2(q^3)\circ \SL_2(q)$ and
the result holds.

In the remaining cases excluding $E_7(q)$, we work in the simply connected group.  Then centralizers of involutions are connected in the algebraic
group. 
If $x$ and $y$ are nonconjugate involutions, then by \cite{Lu},  there is no containment between $C_G(x)$ and $C_G(y)$
and so over the algebraic group, the $C(y)$ orbit of $x$ has positive dimension and so $|x^G \cap C_G(y)|$ grows with $q$.  
If $x$ and $y$ are conjugate, one just notes that there exists a maximal rank semisimple subgroup containing a maximal torus
intersecting $x$ with $x$ noncentral and so the result follows by induction.

In the case of $E_7(q)$, we work in the adjoint group and the connected component of the centralizer of any involution has
index at most $2$ in the centralizer.   The same argument applies working with the connected component of the centralizers. 

\subsection{Exceptional Groups in Even Characteristic}  

Here we work in the adjoint group.   If $G={^2}B_2(q)$, there is a unique class of involutions, the number of conjugates that
it commutes with is $q$.   In $G_2(q)$, there are two classes of involutions corresponding to the two types of root subgroups.
We can choose them to be commuting and so $\|x^G \cap C_G(y)| \ge q$.  If $G={^2}F_4(q)'$,  there are two classes of involutions.
Since the center of a Sylow $2$-subgroup has order at least $q$, the result holds if $x$ is a $2$-central involution (since the center
of a Sylow $2$-subgroup intersects only one class).   Suppose $x$ is a noncentral involution commuting with $y$.   We can choose
$y$ so that $y$ is noncentral in $C_G(x)$ and so $|y^{C_G(x)}|$ grows with $q$ (all involutions have representatives over the field of  $2$
elements since centralizers of involutions
are connected \cite[]{GR}.  

In the remaining cases we observe that there exists a proper semisimple subgroup intersecting all classes of involutions and use
induction.  This was noted already in \cite{GR} aside from the case of $E_7$.  Even in that case it was noted that there are
$4$ commuting simple root subgroups such that the abelian group they generate intersects each class of involutions and
these root subgroups lie in the derived subgroup of a Levi subgroup and so the results holds in that case as well.

\section{Proof of Theorem \ref{t:main2}}

\medskip
Let $n$ be a positive  integer, $\epsilon >0$ be a real number. Suppose that the finite simple group $G$ contains an involution $t$ with $$[C_{G}(t):O_{2^{\prime}}(C_{G}(t))] = n$$ and also that 
$$|G| > (1+\epsilon )|C_{G}(t)|^{2}|\mathcal{I}(C_{G}(t))|.$$ We will prove that there are only finitely many possibilities for $G$.

\medskip
Let us count the number of times the product of two conjugates of $t$ is an element of the form
$tu,$ where $u$ is an element of odd order in $C_{G}(t)$. This is the number of times the product of two conjugates of $t$ lands in ${\rm Sec}_{2}^{C_{G}(t)}(t),$ the $2$-section of $t$ in $C_{G}(t).$ This is the same as the number of $G$-conjugates $t_{1}$ of $t$ and elements $tu \in {\rm Sec}_{2}^{C_{G}(t)}(t)$ such that $t_{1}tu$ is also an involution (of $C_{G}(t)$) which is $G$-conjugate to $t$.

\medskip
This is certainly less than or equal to $$|\mathcal{I}(C_{G}(t))| \times |{\rm Sec}_{2}^{C_{G}(t)}(t)|.$$

\medskip
By the well-known Burnside-Frobenius character-theoretic formula, this number may also be expressed as 

$$\frac{|G|}{|C_{G}(t)|^{2}} \left( \sum_{\chi \in {\rm Irr}(G)}\frac{ \chi(t)^{2}}{\chi(1)}  \sum_{tu \in {\rm Sec}_{2}^{C_{G}(t)}(t)}   \overline{\chi(tu)}  \right).$$

\medskip
Now by Brauer's second main theorem, we have  $$\sum_{tu \in {\rm Sec}_{2}^{C_{G}(t)}(t)}   \overline{\chi(tu)} =0$$ unless $\chi$ lies in $B = B_{0}^{2}(G),$ 
the principal $2$-block of $G.$

\medskip
Hence we may rewrite the expression above as $$\frac{|G|}{|C_{G}(t)|^{2}} \left( \sum_{\chi \in {\rm Irr}(B)}\frac{ \chi(t)^{2}}{\chi(1)}  \sum_{tu \in {\rm Sec}_{2}^{C_{G}(t)}(t)}   \overline{\chi(tu)}  \right).$$

\medskip
Since this expression is bounded above by $$|\mathcal{I}(C_{G}(t)||{\rm Sec}_{2}^{C_{G}(t)}(t)|,$$ there must be some $tu \in {\rm Sec}_{2}^{C_{G}(t)}(t)$ such that 
$$\frac{|G|}{|C_{G}(t)|^{2}} \left( \sum_{\chi \in {\rm Irr}(B)}\frac{ \chi(t)^{2}\overline{\chi(tu)} }{\chi(1)} \right) \leq |\mathcal{I}(C_{G}(t))|.$$ For notice that 
$$\frac{|G|}{|C_{G}(t)|^{2}} \left( \sum_{\chi \in {\rm Irr}(B)}\frac{ \chi(t)^{2}\overline{\chi(tu)} }{\chi(1)}    \right)$$ is always a real number, though, a priori, it might not be non-negative.

\medskip
Hence, for this $2$-regular $u \in C_{G}(t),$ we have $$\frac{|G|}{|C_{G}(t)|^{2}} \left(  1 - \sum_{1 \neq \chi \in {\rm Irr}(B)} \frac{ \chi(t)^{2}|\chi(tu)|}{d} \right) \leq |\mathcal{I}(C_{G}(t))|,$$ where $d$ is the minimal degree of a non-trivial irreducible character in $B$.

\medskip
By Brauer's second and third main theorems, we have $$\sum_{\chi \in {\rm Irr}(B)}\chi(t)^{2} = {\rm dim}(b),$$ where $b$ is the principal $2$-block of $C_{G}(t)$, and this is certainly at most $$[C_{G}(t):O_{2^{\prime}}(C_{G}(t))].$$ Also, (by the same theorems), we have 
$$\sum_{\chi \in {\rm Irr}(B)}|\chi(tu)|^{2} \leq \sum_{\chi \in {\rm Irr}(B)}\chi(t)^{2} \leq [C_{G}(t):O_{2^{\prime}}(C_{G}(t))].$$

In particular, $$|\chi(tu)| <  \sqrt{ [C_{G}(t):O_{2^{\prime}}(C_{G}(t))]}$$ for each $\chi \in {\rm Irr}(B).$

\medskip
Hence we certainly have $$ 1 < \frac{|C_{G}(t)|^{2}|\mathcal{I}(C_{G}(t))|}{|G|} + \frac{[C_{G}(t):O_{2^{\prime}}(C_{G}(t))]^{\frac{3}{2}}}{d}.$$ Since $$|G| >(1+\epsilon)|C_{G}(t)|^{2}|\mathcal{I}(C_{G}(t))|,$$ this yields $$d < (1+\frac{1}{\epsilon})n^{\frac{3}{2}}.$$

\medskip
By Jordan's Theorem on finite subgroups of ${\rm GL}(m,\mathbb{C}),$ there are only finitely many non-Abelian simple groups $G$ which have a non-trivial complex irreducible character of degree less than $$(1+ \frac{1}{\epsilon})n^{\frac{3}{2}},$$ so there are only finitely many possibilities for $G$, as claimed.

\medskip
By contrast, we have:

\medskip
\noindent {\bf Example:} Let $a>1$ be an integer, and let $q$ be a prime power with\\ 
$q \equiv 1+2^{a}$ (mod $2^{a+1}).$ Let $G = {\rm PSL}(2,q).$ Then $G$ has a single conjugacy class of involutions, and if $t$ is an involution of $G$, then the subgroup $C_{G}(t)$ is a dihedral group with $q-1$ elements. Now $q-1$ is divisible by $2^{a},$ but is not divisible by $2^{a+1},$ so we have 
$$[C_{G}(t):O_{2^{\prime}}(C_{G}(t))] = 2^{a}.$$

\medskip
Also, since $q-1$ is divisible by $4$, we see that $C_{G}(t)$ contains $\frac{q+1}{2}$ involutions, so we have $$|G| = \frac{(q-1)q(q+1)}{2} > |C_{G}(t)|^{2}|\mathcal{I}(C_{G}(t))|.$$

 \end{document}